\documentclass[a4paper,12pt]{article} 

\usepackage{amsmath,amssymb,color} 
\usepackage{latexsym} 
\usepackage{graphicx} 

\newtheorem{theorem}{Theorem}

\newtheorem{proposition}[theorem]{Proposition}

\date{\empty}

\begin{document} 

\title{Spanning path-cycle systems with given end-vertices in regular graphs (full version)} 

\author{Yoshimi Egawa$^1$, 
Mikio Kano$^2$
\footnote{mikio.kano.math@vc.ibaraki.ac.jp}
and Kenta Ozeki$^3$
\footnote{ozeki-kenta-xr@ynu.ac.jp}
\and 
$^{1}$  Tokyo University of Sciences, Shinju-ku, Tokyo, Japan
\and
$^{2}$ Ibaraki University, Hitachi, Ibaraki, Japan 
\and
$^3$ Yokohama National University, Yokohama, Japan
} 

\maketitle

\begin{abstract}  We prove the following theorem. Let $r\ge 4$ be an integer, and $G$ be a $K_{1,r}$-free $r$-edge-connected $r$-regular graph. Then, for every set $W$ of even number of vertices of $G$ such that the distance between any two vertices of $W$  in $G$ is at least 3, $G$ has vertex-disjoint paths and cycles $P_1, \ldots, P_m, C_1,  \ldots, C_n$ such that (i) $V(G)=V(P_1) \cup \cdots \cup V(P_m) \cup V(C_1) \cup \cdots \cup V(C_n)$, (ii) each path $P_i$ connects two vertices of $W$, and (iii) the set of the end-vertices of $P_i$'s is equal to $W$.  A similar result for a 3-regular graph is obtained in [Graphs Combin. {\bf 39} (2023) \#85]. However, our proof is widely different from its proof.
\end{abstract}

\section{Introduction} 
We consider simple graphs, which have neither loops nor multiple edges. 
Let $G$ be a graph with vertex set $V(G)$ and edge set $E(G)$.
For a vertex $v$ of a subgraph $H$ of $G$, we denote by $\deg_H(v)$ the degree of $v$ in $H$, and by $N_H(v)$ the neighborhood of $v$ in $H$. Thus $\deg_H(v)=|N_H(v)|$.
 If every vertex of $G$ has degree $r$, then $G$ is called an {\em $r$-regular graph}.
The {\em star} $K_{1,m}$ of order $m+1$ is the complete bipartite graph with partite sets of size 1 and $m$.
Furthermore, the star $K_{1,3}$ is often called a {\em claw}.
A graph that has no induced subgraph isomorphic to $K_{1,m}$ is called a {\em $K_{1,m}$-free graph}.

Let $W$ be a set of even number of vertices of a graph $G$. Then we say that $G$ has a {\em path system with respect to $W$} if there are vertex-disjoint paths $P_1,P_2, \ldots, P_ m$ in $G$ such that (i) each path $P_i$ connects two vertices of $W$ and (ii) the set of the end-vertices of  $P_i$'s is equal to $W$; in particular, no internal vertex of each $P_i$ is contained in $W$, and $m=|W|/2$. 
Moreover, we say that $G$ has a {\em spanning path-cycle system with respect to $W$} 
if there are vertex-disjoint paths and cycles  $P_1,P_2, \ldots, P_ m, C_1, C_2, \ldots C_n$ in $G$ such that (i) each path $P_i$ connects two vertices of $W$, (ii) the set of end-vertices of  $P_i$'s is equal to $W$, 
(iii) each $C_j$ is a cycle of $G$, and (iv) $V(G)=V(P_1) \cup \cdots \cup V(P_m) \cup V(C_1) \cup \cdots \cup V(C_n)$. In other words, if $\{P_1,P_2, \ldots, P_ m, C_1, C_2, \ldots C_n\}$ is a spanning path-cycle system with respect to $W$, then $\{P_1,P_2, \ldots, P_ m\}$ is a path-system with respect to $W$ and $V(G)- \big( V(P_1) \cup \cdots \cup V(P_m)\big) $ is covered by vertex-disjoint cycles $C_1, C_2, \ldots, C_n$. 
 
 In this paper, we present some results about spanning path-cycle systems. Before giving them, let us begin with some results about path systems.
 
\begin{theorem} [Kaiser \cite{Kaiser-2008}]
Let $r \geq 2$ be an integer, and $G$ be an $r$-edge-connected $r$-regular graph. 
Then, for every set $W$ of even number of vertices of $G$, 
$G$ has a  path system with respect to $W$.
\label{thm-10} 
\end{theorem}

\begin{theorem} [Furuya and Kano \cite{FK-2023}]
Let $G$ be a connected claw-free graph. Then, for every set $W$ of even number of vertices of $G$, 
$G$ has a path system with respect to $W$.
\label{thm-11} 
\end{theorem} 

A criterion for a graph to have a  path-system with respect to $W$ for every set $W$ of even number of vertices is given in the following theorem, and the above two theorems can be proved by using this theorem.
Note that $\omega(G)$ denotes the number of components of $G$.

\begin{theorem}[Lu and Kano \cite{LK-2020}]
Let $G$ be a connected  graph. Then $G$ has a path system with respect to $W$ 
for every set $W$ of even number of vertices of $G$ if and only if
\[ \omega(G-S) \le |S|+1 \quad \mbox{for all} \quad S\subset V(G).\]
\label{thm-12} 
\end{theorem} 

\vspace{-9mm}
We now give some results about spanning path-cycle systems including our theorems.

\begin{theorem} [Furuya and Kano \cite{FK-2023}]
Let $G$ be a claw-free $3$-edge-connected $3$-regular graph. Then, for every set $W$ of even number of vertices such that the distance between any two vertices of  $W$ in $G$ is at least $3$, $G$ has a spanning
path-cycle system with respect to $W$.
\label{thm-1} 
\end{theorem} 

The following is our result.

\begin{theorem} 
Let $r\ge 4$ be an integer, and $G$ be a $K_{1,r}$-free $r$-edge-connected $r$-regular graph. 
Then, for every set $W$ of even number of vertices such that the distance between any two vertices of $W$ in $G$ is at least $3$, $G$ has a spanning path-cycle system with respect to $W$.
\label{thm-2} 
\end{theorem} 

Actually, we prove Theorem~\ref{thm-3}, 
and Theorem~\ref{thm-2} is an easy consequence of Theorem~\ref{thm-3} since if $W$ satisfies the condition in Theorem~\ref{thm-2}, then it also satisfies the condition in Theorem~\ref{thm-3}. 
However, Theorem~\ref{thm-2} is not equivalent to Theorem~\ref{thm-3} because some two vertices of $W$ given in Theorem~\ref{thm-3}
may be joined by edges of $G$.
 
\begin{theorem} 
Let $r\ge 4$ be an integer, and $G$ be a $K_{1,r}$-free $r$-edge-connected $r$-regular graph.
Then, for every set $W$ of even number of vertices such that $|N_G(v)\cap W|\le 1$ for every vertex $v$ of $G$,  $G$ has a spanning path-cycle system with respect to $W$.
\label{thm-3} 
\end{theorem} 

We prove Theorem~\ref{thm-3} in the next section.
In Section 3, we show that  the conditions on $W,$ and the edge-connectivity in Theorems~\ref{thm-2} and \ref{thm-3} are sharp, and also that  the condition of $K_{1,r}$-freeness is necessary.

\section{Proof of Theorem~\ref{thm-3}} 

We begin with some notation and definitions.
Let $G$ be a graph.
An edge joining a vertex $x$ to a vertex $y$ is denoted by $xy$ or $yx$.
For two disjoint vertex sets $X$ and $Y$ of $G$, let
$E_G(X,Y)$ denote the set of edges of $G$ joining $X$ to $Y$, and $e_G(X,Y)$ denote the number of edges of $G$ joining $X$ to $Y$. Thus $e_G(X,Y)=|E_G(X,Y)|$. 
If $X=V(D)$ and $Y\subset V(G)-V(D)$ for some subgraph $D$ of $G$, then we briefly write
$e_G(D,Y)$ for $e_G(V(D),Y)$. 
For a vertex set $X$ of $G$, the subgraph of $G$ induced by $X$ is denoted by $G[X]$.
If a vertex $v$ of $G$ is contained in a subgraph $D$ of $G$, 
then we briefly write $v\in D$ instead of $v\in V(D)$.

In order to prove Theorem~\ref{thm-3},  we focus on an $f$-factor of a graph $G$.
Let $\mathbb{Z}^+$ denote the set of non-negative integers.
For a function $f:V(G) \to \mathbb{Z}^+$, a spanning subgraph $F$ of $G$ is called an {\em $f$-factor} of $G$ if
$ \deg_F(v) = f(v)$ for all $v\in V(G)$.

It is easy to see that a graph $G$ has a spanning path-cycle system with respect to $W$ if and only if  $G$ has a factor $F$ that satisfies  
\begin{align*}
 \deg_F(x) & =1 \quad  \mbox{for every} \quad x\in W, ~\mbox{and} \\
 \deg_F(y) & =2 \quad  \mbox{for every} \quad y\in V(G)-W. 
 \end{align*}
Namely, each component of $F$ is a path or a cycle, and the set of paths and cycles of $F$ forms a spanning path-cycle system with respect to $W$. Thus we prove Theorem~\ref{thm-3} by using the following $f$-factor theorem.

For an integer-valued function $h$ defined on $V(G)$ and a subset $X\subseteq V(G)$, we briefly write
\[ h(X):=\sum_{x\in X} h(x) \quad \mbox{and} \quad \deg_G(X):= \sum_{x\in X} \deg_G(x).\]
A criterion for a graph to have an $f$-factor is given in the following theorem,
which is called ``The $f$-factor Theorem''.

\begin{theorem}[Tutte \cite{Tutte-1952}, \cite{Tutte-1954}, Theorem 3.2 in \cite{AK-2011}] Let $G$ be a graph, and $f:V(G) \to \mathbb{Z}^+$ be a function. Then $G$ has an $f$-factor if and only if for all disjoint subsets $S,T \subseteq V(G)$, 
\begin{align} 
\delta(S,T) := f(S) + \deg_{G-S}(T) -f(T) -q(S,T) \ge 0,
 \label{eq-2} 
\end{align} 
where $q(S,T)$ denotes the number of components $D$ of $G-(S\cup T)$ satisfying 
\begin{align}  f(V(D))+ e_G(V(D),T) \equiv 1 \pmod{2}. 
\label{eq-3} 
\end{align} 
\label{thm-10} 
In addition,
we have
$\delta(S,T) \equiv f(V(G)) \pmod{2}$.
\end{theorem}

 \noindent 
Note that a component $D$ of $G-(S\cup T)$ satisfying (\ref{eq-3}) is called an {\em $f$-odd component} of $G-(S \cup T)$.  

\bigskip \noindent
{\em Proof of Theorem~\ref{thm-3}.}  Define a function $f :V(G) \to \mathbb{Z}^+$ by letting 
\[ f(v)= \bigg\{
\begin{array}{ll} 
1 & \mbox{if ~$v \in W$,} \\ 
2 & \mbox{if ~$v \in V(G)-W$.}
\end{array} 
    \] 
Then $G$ has the desired spanning path-cycle system with respect to $W$  if and only if $G$ has an $f$-factor.

Assume that $G$ has no $f$-factor.
Then, by Theorem \ref{thm-10}, there exist two disjoint vertex sets $S$ and $T$ of $G$ such that
$\delta(S,T)=f(S) + \deg_{G-S}(T) -f(T)  -q(S,T) < 0$.
We take such $S$ and $T$ so that $|T|$ is as small as possible.

Since $f(V(G))$ is even, $\delta(\emptyset, \emptyset)= -q(\emptyset, \emptyset)=0$, which implies that $S\cup T \ne \emptyset$.

\medskip \noindent
{\bf Claim~1.} {\em $S \ne \emptyset$ and $T \ne \emptyset$.}

\medskip \noindent
{\em Proof.}  Assume that $T=\emptyset$. Then $S\ne \emptyset$ as $S\cup T \ne\emptyset$. Let $D_1,D_2, \ldots, D_m$ be the $f$-odd components of $G-S$, where $m= q(S, \emptyset)$. Then $e_G(S, D_i) \ge r$ for every $1\le i \le m$ by the edge connectivity of $G$,  and so $r|S| =\sum_{x\in S}\deg_G(x) \ge \sum_{1\le i \le m} e_G(S, D_i) \ge rm$, which implies $|S|\ge m$. Thus $\delta(S,\emptyset)= f(S) -q(S, \emptyset ) \ge |S| -m \ge 0$, which contradicts our choice of $S$ and $T$. Thus $T\ne \emptyset$.

Next assume that $S=\emptyset$. Then $T\ne \emptyset$ as $S\cup T \ne\emptyset$. Let $D_1,D_2, \ldots, D_{m'}$ be the $f$-odd components of $G-T$, where $m'=q(\emptyset, T)$. By the same argument as given above, we have $|T|\ge m'$. Then $\delta(\emptyset, T)= \deg_G(T) -f(T) -q(\emptyset, T) \ge r|T| -2|T| -m' \ge  0$, a contradiction.  Thus $S\ne \emptyset$.~~$\Box$
\\

\medskip \noindent
{\bf Claim~2.} {\em No two vertices in $T$ are adjacent in $G$.}

\medskip \noindent
{\em Proof.}  Assume that two vertices $y, y' \in T$ are adjacent in $G$.
Let $T' = T- \{y'\}$.
Note that 
$\deg_{G-S}(T') = \deg_{G-S}(T) - \deg_{G-S}(y')$,
$f(T') = f(T) - f(y')$,
and $q(S,T') \geq q(S,T) - q'$,
where $q'$ is the number of $f$-odd components $D$ of $G-(S\cup T)$
such that $D$ is adjacent to $y'$.
Since $y'$ is adjacent to $y \in T$,
we have $\deg_{G-S}(y') \geq q' +1$.
Since $f(y') \leq 2$, we  therefore obtain
\begin{align*} 
\delta(S,T') &= 
f(S) + \deg_{G-S}(T') -f(T') -q(S,T') \\
&\leq f(S) + \deg_{G-S}(T) - \deg_{G-S}(y')  \\
& \hspace{2em} - \big( f(T) - f(y') \big) - \big( q(S,T)-q' \big) \\
&\leq \delta(S,T) + f(y') -1 \\
&\leq \delta(S,T)  + 1.
\end{align*} 
By Theorem \ref{thm-10},
we see 
$\delta(S,T) \equiv f(V(G)) \equiv \delta(S,T') \pmod{2}$,
which implies that 
$\delta(S,T') \leq \delta(S,T) < 0$.
However, this contradicts the minimality of $T$.~~$\Box$
\\

\medskip 
Let $S_i = \{x \in S : f(x) = i\}$ for $i = 1,2$,  and $\mathcal{D}=\{D_1,\ldots, D_m\}$ be the set of $f$-odd components of $G-(S\cup T)$, where $m=q(S,T)$, 
and let $U=V(D_1)\cup \cdots \cup V(D_m)$.
We use a discharging method to prove Theorem~\ref{thm-3}.
We set a function $\varphi : S \cup T \cup \mathcal{D} \to \mathbb{R}$, as an initial charge,
as follows:
For every  $D\in \mathcal{D}$, let $\varphi(D)  = 0$, and for each $v \in S \cup T$,
let
\[
 \varphi(v) = 
\begin{cases}
1 & \text{if $v \in S_1$}, \\
2 & \text{if $v \in S_2$}, \\
\deg_{G-S-U}(v) & \text{if $v \in T$}.
\end{cases} 
\]
Note that 
$\deg_{G-S}(T) = \deg_{G-S-U}(T) + e_G(T,U)$,
and hence
\begin{align}
\sum_{x \in S \cup T} \varphi(x) &= f(S) + \deg_{G-S-U}(T) 
\nonumber\\
&= f(S) + \deg_{G-S}(T) - e_G(T,U).
 \label{eq-charge}
\end{align}
We now design some discharging rules to redistribute charges
between vertices in $S \cup T$ and $\mathcal{D}$ along  edges as follows.

\begin{itemize}
\item[(i)]
Let $xy$ be an edge joining $x \in S_1$ to $y \in T \cup U$.
Then,
$x$ sends a charge of $1/r$ to $y$ when $y \in T$, 
and to the $f$-odd component containing $y$ when $y \in U$.
\item[(ii)]
Let $xy$ be an edge joining $x \in S_2$ to $y \in T \cup U$.
\begin{itemize}
\item[(ii-1)]
If $y \in T$,
then $x$ sends a charge of  $\displaystyle \frac{2r-1}{r(r-1)}$ to $y$.
\item[(ii-2)]
If $y \in U$, 
then $x$ sends a charge of $1/r$ to the $f$-odd component containing $y$.
\end{itemize}
\item[(iii)]
Let $zy$ be an edge joining $z \in U$ to $y \in T$.
Then the $f$-odd component containing $z$ sends a charge of  $\displaystyle \frac{r-1}{r}$ to $y$.
\end{itemize}

\noindent
It is easy to see that the following inequalities hold because $r\ge 4$.
\begin{align}
\frac{1}{r} \leq 
\frac{2r-1}{r(r-1)} \leq \frac{r-1}{r}.
\label{basic_ineq}
\end{align}

For each $v \in S \cup T$ 
and for each $D\in \mathcal{D}$,
let $\varphi^*(v)$ and $\varphi^*(D)$, respectively, denote the charge of $v$ and $D$
after the discharging procedure.
Then the following claims hold.

\medskip \noindent
{\bf Claim~3.} {\em For each $x \in S$, we have $\varphi^*(x) \geq 0$.}

\medskip  \noindent 
{\em Proof.}  Let $x \in S$.
If $x \in S_1$, then by rule (i),
$x$ sends a charge of $1/r$ to each of its neighbors
belonging to $T \cup \mathcal{D}$. Thus,
\[ \varphi^*(x) \geq \varphi(x) - \frac{1}{r} \times r = 1-1 = 0. \]

Suppose next $x \in S_2$. Then $\varphi(x) = 2$.
Since $G$ is $K_{1,r}$-free,
some two vertices in the neighborhood of $x$ must be adjacent.
Thus, 
if all neighbors of $x$ belong to $T$,
then this contradicts Claim~2.
Therefore, $x$ has a neighbor $y_1 \notin T$.
In this case,
$x$ sends a charge of  at most $1/r$ to $y_1$,
and hence by (ii) and (\ref{basic_ineq}), 
$$\varphi^*(x) \geq 2- \frac{1}{r} - \frac{2r-1}{r(r-1)} \times (r-1) = 0,$$
as desired.
Hence the claim holds. ~~$\Box$

\medskip \noindent 
{\bf Claim~4.} {\em For each $y \in T$, we have $\varphi^*(y) \geq 2$.}

\medskip  \noindent
{\em Proof.} Let $y \in T$.
By the rules,
$y$ receives only positive charges,
and hence $\varphi^*(y) \geq \varphi(y)$.
If $\deg_{G-S-U}(y) \geq 2$,
then 
$\varphi^*(y) \geq \varphi(y) = \deg_{G-S-U}(y) \geq 2$,
and we are done.
Thus, we may assume $\deg_{G-S-U}(y) \leq 1$.
This implies that $y$ is adjacent to at least $r-1$ vertices in $S \cup U$.
In addition, $N_G(y)$ contains at most one vertex of $S_1 \subseteq W$ by the assumption of the theorem.
Hence,  by (\ref{basic_ineq}), the worst case is the case where one neighbor of $y$ is in $S_1$ and all the other neighbors of $y$ in $S\cup T$ lie in $S_2$.
Therefore, if $\deg_{G-S-U}(y) = 1$, then 
\begin{align*}
\varphi^*(y) &\geq \deg_{G-S-U}(y) + \frac{2r-1}{r(r-1)} \times (r-2) + \frac{1}{r} 
\\
&= \frac{3r^2-5r+1}{r(r-1)}
= 2 + \frac{(r-\frac{3}{2})^2 - \frac{5}{4}}{r(r-1)}
\geq 2.
\end{align*}
If $\deg_{G-S-U}(y) = 0$, then 
$$\varphi^*(y) \geq \frac{2r-1}{r(r-1)} \times (r-1) + \frac{1}{r} = 2,$$
and we are done. ~~$\Box$

\medskip \noindent
{\bf Claim~5.} {\em  For each $D\in \mathcal{D}$,
we have $\varphi^*(D) \geq 1 - e_G(T,D)$.}

\medskip  \noindent 
{\em Proof.} Let $D \in \mathcal{D}$,
and let $x_1z_1, \dots , x_az_a, y_{1}z_{a+1}, \dots , y_{b}z_{a+b}$
be the edges joining $S \cup T$ to $D$,
where $x_1, \ldots, x_a \in S$, $y_{1} ,\ldots, y_{b}\in T$ and $z_1, \ldots, z_{a+b}\in D$.
Note that $a=e_G(S,D)$ and $b = e_G(T,D)$.
Since $G$ is $r$-edge-connected,
we have $a+b \geq r$.
By the rules (i) and (ii-2), every $x_i$ sends a charge of $1/r$ to $D$,
and  by the rule (iii), $D$ sends a charge of $(r-1)/r$ to every $y_j$.
Thus, we have 
\[ \varphi^*(D) = \frac{1}{r} \times a - \frac{r-1}{r} \times b
= \frac{a+b}{r} - b
\geq 1 - e_G(T,D). \]
Therefore the claim follows. ~~$\Box$

\bigskip

By Claims~3 and 4,
we have 
$\sum_{x \in S \cup T} \varphi^*(x) \geq 2|T| \geq f(T)$.
By Claim~5,
we have 
$\sum_{D \in \mathcal{D}} \varphi^*(D) \geq  q(S,T) - e_G(T,U)$.
Since the sum of charges is preserved by the discharging step,
it follows from (\ref{eq-charge}) and the above two inequalities 
that
\begin{align*}
\delta(S,T) & = f(S) + \deg_{G-S}(T) -f(T)  -q(S,T) \\
& = \sum_{x \in S \cup T} \varphi(x) +  e_{G}(T,U) -f(T) -q(S,T) \\
& = \sum_{x \in S \cup T} \varphi^*(x) + \sum_{D \in \mathcal{D}} \varphi^*(D)  
  + e_{G}(T,U) -f(T) -q(S,T)\\
& \ge 0.
\end{align*}
Consequently, Theorem~\ref{thm-3} is proved. ~~$\Box$

\section{Sharpness of Theorems~\ref{thm-2} and \ref{thm-3}}

In this section, we show that some conditions in Theorems~\ref{thm-2} and \ref{thm-3} are sharp or necessary.  Namely, we show that  (i) $r$-edge-connectedness cannot be replaced by $(r-1)$-edge-connectedness when $r$ is odd,
and by $(r-2)$-edge-connectedness when $r$ is even,
 (ii) the condition that $G$ is $K_{1,r}$-free cannot be removed, and (iii) the condition on $W$ cannot be replaced by a weaker condition.
Note that the sharpness of Theorem~\ref{thm-1} is shown in \cite{FK-2023}. 

We here note that if $r$ is even, then every $(r-1)$-edge-connected $r$-regular graph is $r$-edge-connected. Thus, if $r$ is even, then in order to show the sharpness of $r$-edge-connectivity, it suffices to verify that the desired conclusion does not hold for $(r-2)$-edge-connected $r$-regular graphs. 

We first prove the following proposition.

\begin{proposition} \label{prop-1}
Let $r\ge 4$ be an integer. Then the following statements hold, where $W$ denotes a set of even number of vertices of $G$ such that the distance between any two vertices of $W$ in $G$ is at least 3.
\begin{itemize}
\item[(1)] If $r$ is odd, then there are infinitely many pairs $(G,W)$ of a  $K_{1,r}$-free $(r-1)$-edge-connected $r$-regular graph $G$ and $W\subset V(G)$ such that $G$ has no spanning path-cycle system with respect to $W.$

\item[(2)]  If $r$ is even, then there are infinitely many pairs $(G,W)$ of a $K_{1,r}$-free $(r-2)$-edge-connected $r$-regular graph $G$ and $W\subset V(G)$ such that $G$ has no spanning path-cycle system with respect to $W.$

\item[(3)]  For every $r$, there are infinitely many pairs $(G,W)$ of an $r$-edge-connected $r$-regular graph $G$ and $W\subset V(G)$ such that $G$ has no spanning path-cycle system with respect to $W.$
\end{itemize}
\end{proposition}

\begin{figure}[htbp]
\begin{center}
\includegraphics[scale=0.8]{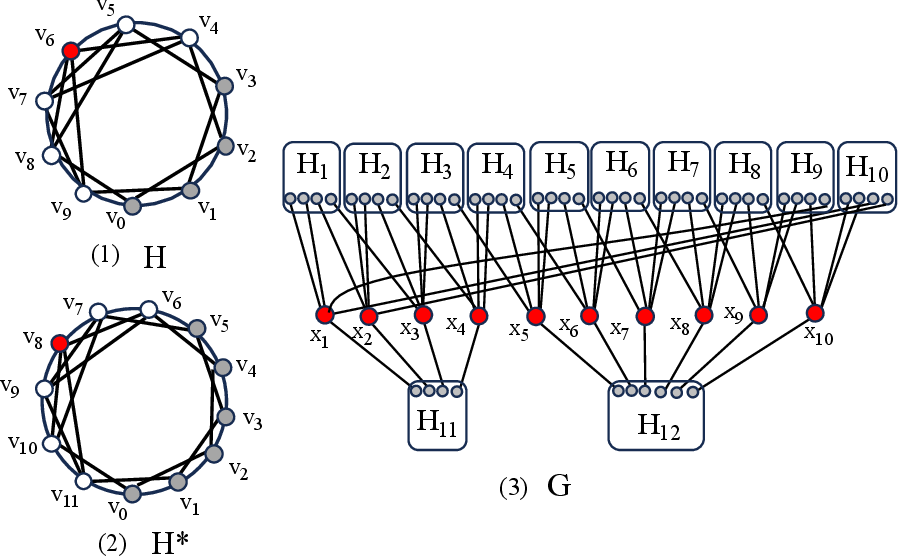}
\caption{(1) The graph $H$ with $r=5$ and $k=6$, where $v_6=v_{r+1}$, every grey vertex has degree $r-1=4$, and all the other vertices have degree $r=5$.  (2) The graph $H^*$ with $r=5$ and $k=6$, where $v_8=v_{r+3}$, every grey vertex has degree $r-1=4$, and all the other vertices have degree $r=5$. 
(3) A $K_{1,r}$-free $(r-1)$-edge-connected $r$-regular graph $G$ that has no spanning path-cycle system with respect to $W=\{x_1,x_2, \ldots x_{10}\} \cup \{v_6\in H_i :1\le i \le 11\} \cup
 \{v_{8} \in H_{12} \}$.}
\label{fig-2}
\end{center}
\end{figure}

\noindent 
{\em Proof.} We first prove (1). Let $r\ge 5$ be an odd integer, and $k\ge r+1$ be an even integer. We define a graph $H$ with vertex set $\{v_0,v_1, \ldots, v_{r+k-2}\}$ as follows. For convenience, let $v_i=v_j$ if $i\equiv j \pmod{r+k-1}$. The  edge set of $H$ is 
\begin{align*}
 E(H) = & ~ \{ v_iv_j : |i-j| \le (r-1)/2 \}  \\
 &  \cup \{v_{s}v_{s+ k/2} : r-1\le s \le r+k/2-2\}. \hspace{2em} \mbox{ (see (1) of Fig.~\ref{fig-2})} 
 \end{align*} 
Then $H$ has $r-1$ vertices $v_0,v_1, \ldots, v_{r-2}$ with degree $r-1$
and $k$ vertices $v_{r-1}, v_{r}, \ldots, v_{r+k-2}$ with degree $r$. 
 A graph $H^*$ with vertex set $\{v_0,v_1, \ldots, v_{r+k}\}$ is defined as follows.
\begin{align*}
 E(H^*) = & ~ \{ v_iv_j : |i-j| \le (r-1)/2 \}   \\
  &  \cup \{v_{s}v_{s+k/2} : r+1\le s \le r+k/2\}. \hspace{2em} \mbox{ (see (2) of Fig.~\ref{fig-2})}  
 \end{align*}
 Then $H^*$ has $r+1$ vertices $v_0,v_1, \ldots, v_{r}$ with degree $r-1$
and $k$ vertices $v_{r+1}, v_{r+2}, \ldots, v_{r+k}$ with degree $r$.

Let $H_1,H_2, \ldots, H_{2r+1}$ be $2r+1$ disjoint copies of $H$, and let $H_{2r+2} =H^*$.
We now construct the desired $K_{1,r}$-free $(r-1)$-edge-connected $r$-regular graph $G$.
Let $V(G)=\{x_1,x_2, \ldots, x_{2r}\}\cup V(H_1) \cup V(H_2) \cup \cdots \cup V(H_{2r+2}) $.
For every $H_i, 1\le i \le 2r$, add $r-1$ edges  $v_0x_i$, $v_1x_i$, $v_2x_{i+1}$, $v_3x_{i+2}$, $\ldots$, $v_{r-2}x_{i+r-3}$, where $v_0, v_1, \ldots, v_{r-2}$ are as in the definition of $H$ and the indies of $x$ are taken modulo $2r$.
Then join every vertex of $H_{2r+1}\cup H_{2r+2}$ with degree $r-1$ to a vertex in $\{x_1,x_2, \ldots, x_{2r}\}$ so that the resulting graph $G$ becomes an $r$-regular graph (see (3) of Fig.~\ref{fig-2}).
Then $G$ is  a $K_{1,r}$-free $(r-1)$-edge-connected $r$-regular graph.

Let $W =\{x_1, x_2, \ldots, x_{2r}\} \cup \{ v_{r+\frac{k}{2} -2}\in V(H_j): 1\le j \le 2r+1\} \cup \{v_{r+\frac{k}{2}} \in V(H_{2r+2})\}$. Then the distance between any two vertices of $W$ is at least 3. Moreover, $G$ has no spanning path-cycle system with respect to $W$. To see this, apply Theorem 7 with $S=\{x_1,x_2, \ldots, x_{2r}\}$ and $T=\emptyset$ and with $f$ as in the proof  of Theorem 6. Then $f(S)=2r$ and $q(S,T)=2r+2$. Hence
$\delta(S,T)=-2$, which implies that there is no such system.

\begin{figure}[htbp]
\begin{center}
\includegraphics[scale=0.8]{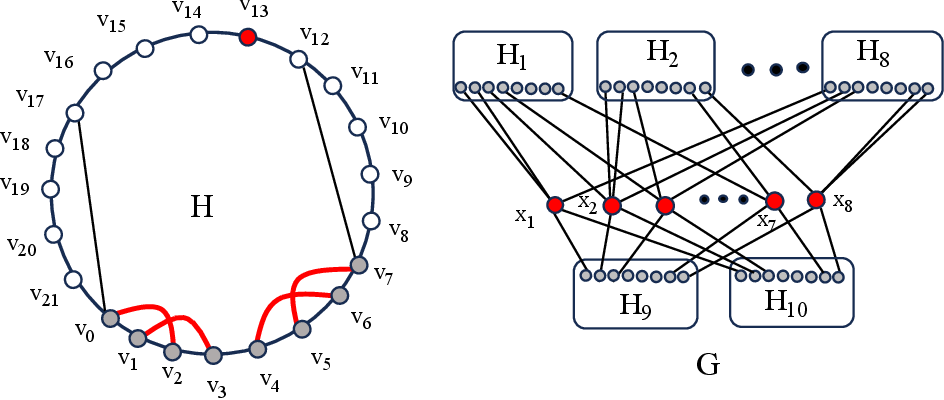} 
\caption{The graph $H$ with $r=10$ and $k=12$,  where $v_0v_2, v_1v_3, v_4v_6,  v_5v_7, \not\in E(H)$, $v_0v_{17}, v_{7}v_{12}\in E(H)$ and $w=v_{13}$. A $K_{1,r}$-free $(r-2)$-edge-connected $r$-regular graph $G$ that has no spanning path-cycle system with respect to 
$W=\{x_1,x_2, \ldots x_{r-2}\} \cup \{w\in H_i :1\le i \le r \} $.}
\label{fig-2C}
\end{center}
\end{figure}

We next prove (2).   We first consider the case where $r=4m+2 \ge 6$ and $k\ge r$.
We define a graph $H$ as follows: $V(H)=\{v_0,v_1, \ldots, v_{r+k}\}$ and 
\begin{align*}
 E(H) &  =  ~ \{ v_iv_j : |i-j| \le r/2 \}   \\
 &\quad  -  \{v_{0}v_2, v_1v_3, v_4v_6, v_5v_7, \ldots, v_{r-6}v_{r-4}, v_{r-5}v_{r-3} \}  .
 \end{align*}   
Then $H$ has $r-2$ vertices $v_0,v_1, \ldots, v_{r-5}, v_{r-4}, v_{r-3}$ with degree $r-1$, and all the other vertices have degree $r$. Put $w=v_{(3r-4)/2}$, which is not adjacent to $v_0,v_1, \ldots, v_{r-3}$  (see Fig.~\ref{fig-2C}).

Let $H_1,H_2, \ldots, H_{r}$ be $r$ disjoint copies of $H$.
We construct a $K_{1,r}$-free $(r-2)$-edge-connected $r$-regular graph $G$ as follows:
Let $V(G)= \{x_1,x_2, \ldots, x_{r-2}\} \cup V(H_1) \cup V(H_2) \cup \cdots \cup V(H_{r}) $.
For each $H_i, 1\le i \le r-2$, add $r-2$ edges  
$v_0x_i$, $v_1x_i$, $v_2x_{i+1}$, $v_3x_{i+2}$, $\ldots$, $v_{r-3}x_{i+r-4}$, where $v_0, v_1, \ldots, v_{r-3}$ are as in the definition of $H$ and the indies of $x$ are taken modulo $r-2$.
Additionally, for every $H_j \in \{H_{r-1},H_{r}\}$, add $r-2$ edges  
$v_0x_1$, $v_1x_2, \ldots, v_{r-3}x_{r-2}$ (Fig.~\ref{fig-2C}).

Let $W =\{x_1, x_2, \ldots, x_{r-2}\} \cup \{ w\in V(H_j): 1\le j \le r\}$. Then the distance between any two vertices of $W$ is at least 3. Moreover, $G$ has no spanning path-cycle system with respect to $W$. To see this, apply Theorem 7 with $S=\{x_1, \ldots, x_{r-2}\}$ and $T=\emptyset$ and with $f$ as in the proof  of Theorem 6. Then $f(S)=r-2$ and $q(S,T)=r$. Hence $\delta(S,T)=-2$, which implies that there is no such system.

\begin{figure}[htbp]
\begin{center}
\includegraphics[scale=0.8]{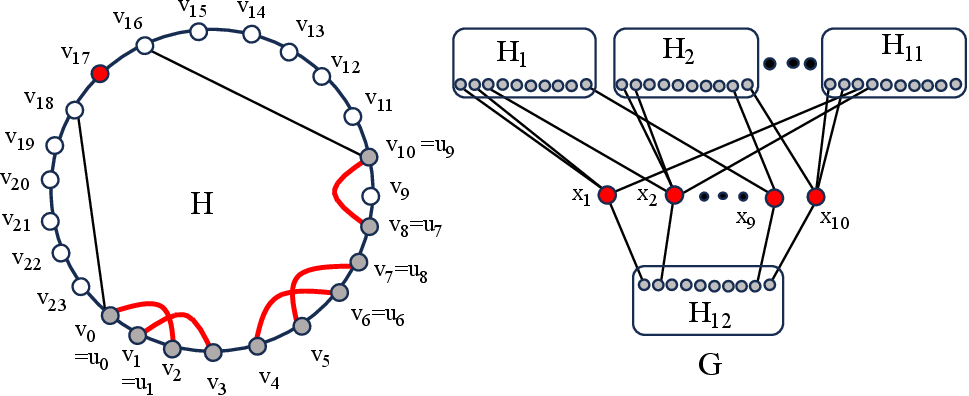}
\caption{The graph $H$ with $r=12$ and $k=12$,  where $v_0v_2, v_1v_3, v_4v_6,  v_5v_7, v_8v_{10} \not\in E(H)$, $v_0v_{18}, v_{10}v_{16}\in E(H)$ and $w=v_{17}$. A $K_{1,r}$-free $(r-2)$-edge-connected $r$-regular graph $G$ that has no spanning path-cycle system with respect to 
$W=\{x_1,x_2, \ldots x_{r-2}\} \cup \{w\in H_i :1\le i \le r \} $.}
\label{fig-2B}
\end{center}
\end{figure}

 We next consider the case where $r=4m\ge 4$. Let $k\ge r$ be an even integer. We define a graph $H$ with vertex set $\{v_0,v_1, \ldots, v_{r+k-1}\}$ and edge set
\begin{align*}
 E(H) &  =  ~ \{ v_iv_j : |i-j| \le r/2 \}   \\
 &\quad  - \Big( \{v_0v_2, v_1v_3, v_4v_6, v_5v_7, \ldots, v_{r-8}v_{r-6}, v_{r-7}v_{r-5}\} \cup \{v_{r-4}v_{r-2} \} \Big) .
 \end{align*}   
Then $H$ has $r-2$ vertices $v_0,v_1, \ldots,  v_{r-4}, v_{r-2}$ with degree $r-1$,
and all the other vertices have degree $r$. 
Put $w=v_{(3r-2)/2}$, which is not adjacent to $v_0,v_1, \ldots, v_{r-2}$ (see Fig.~\ref{fig-2B}). 

Let $H_1,H_2, \ldots, H_{r}$ be $r$ disjoint copies of $H$.
We arrange the vertices of $H_i$ with degree $r-2$ as follows so that every two consecutive vertices $u_j, u_{j+1}, 0\le j \le r-4$,  are adjacent in $H_i$.
\begin{align}
 u_0=v_0, u_1=v_1, \ldots   , u_{r-6}=v_{r-6},  u_{r-5}=v_{r-4},   u_{r-4}=v_{r-5}, u_{r-3}=v_{r-2}.
 \label{eq-11}
\end{align}
We construct a $K_{1,r}$-free $(r-2)$-edge-connected $r$-regular graph $G$ as follows:
Let $V(G)= \{x_1,x_2, \ldots, x_{r-2}\} \cup V(H_1) \cup V(H_2) \cup \cdots \cup V(H_{r}) $.
For each $H_i, 1\le i \le r-1$, add $r-1$ edges  
\[u_0x_i, u_1x_i, u_2x_{i+1}, u_3x_{i+2}, \ldots, 
u_{r-4}x_{i+r-5}, u_{r-3}x_{i+r-4},\] 
where $u_0, u_1, \ldots u_{r-3}$ are defined in (\ref{eq-11}) and the indies of $x$ are taken modulo $r-2$.
Moreover, for $H_{r}$, add $r-2$ edges  
$u_0x_1$, $u_1x_2, \ldots, u_{r-3}x_{r-2}$ (see Fig.~\ref{fig-2B}).

Let $W =\{x_1, x_2, \ldots, x_{r-2}\} \cup \{ w\in V(H_j): 1\le j \le r\}$. Then the distance between any two vertices of $W$ is at least 3. Moreover, $G$ has no spanning path-cycle system with respect to $W$. To see this, apply Theorem 7 with $S=\{x_1, \ldots, x_{r-2}\}$ and $T=\emptyset$ and with $f$ as in the proof  of Theorem 6. Then $f(S)=r-2$ and $q(S,T)=r$. Hence $\delta(S,T)=-2$, which implies that there is no such system.

Now we prove (3). Let $G$ be an $r$-edge-connected $r$-regular bipartite graph with
sufficiently large order. Then $G$ is not $K_{1,r}$-free and has two vertices with distance 4 in the same partite set, and let $W$ be the set of these two vertices. Then, since the two partite sets of $G$ have the same size, $G$ has no spanning path-cycle system with respect to $W$. Hence (3) is proved. ~~$\Box$

\medskip

The following proposition shows that the condition on $W$ in Theorem~\ref{thm-3}  is sharp.

\begin{proposition} \label{prop-2}
Let $r\ge 4$ be an integer. Then there are infinitely many $K_{1,r}$-free $r$-edge-connected $r$-regular graphs $G$ such that there exists a set $W$ of even number of independent vertices in $G$ satisfying  $|N_G(v)\cap W| \le 2$ for every $v\in V(G)$, for which $G$ has no spanning path-cycle system with respect to $W$.
\end{proposition}

\noindent 
{\em Proof of Proposition~\ref{prop-2}.}
A graph is said to be 
\emph{essentially $k$-edge-connected}
if removing any at most $k-1$ edges
results in a graph having at most one component of order at least two.

\begin{figure}[htbp]
\begin{center}
\includegraphics[scale=0.8]{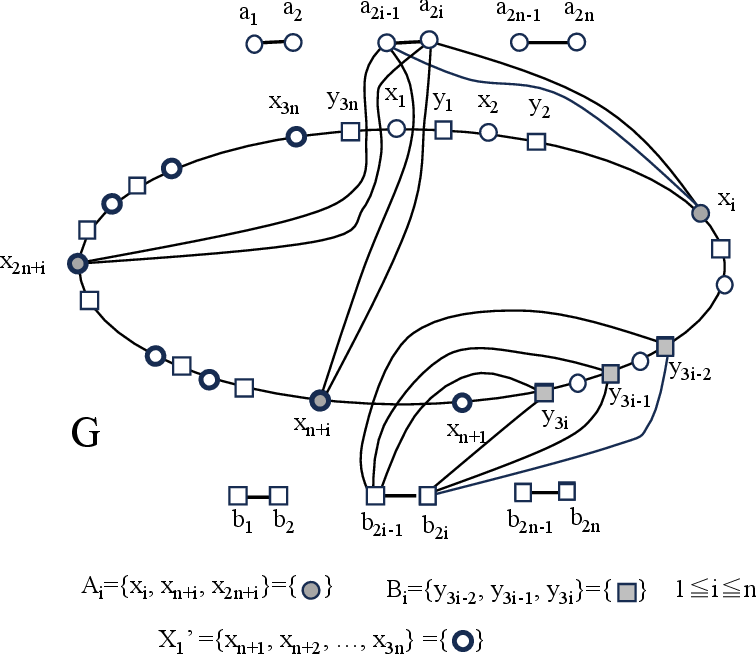}
\caption{A $K_{1,4}$-free 4-edge-connected 4-regular graph $G$. }
\label{fig-6}
\end{center}
\end{figure}

We first consider the case where $r\ge 6$, and later, we deal with the case where $r=4,5$.
Note that if $r=4$, then following bipartite graph $H_2$ cannot be defined, and if $r=5$,
then the desired graph becomes a graph shown in Fig.~\ref{fig-6},
and we need some improved argument given in the proof of the case where $r=5$.

Let $m$ be a multiple of $r-1$
such that $m \geq 2(r-1)^2$.
We first prepare two bipartite graphs as follows:
\begin{itemize}
\item
Let $H_1$ be an 
essentially $3$-edge-connected 
bipartite graph 
with bipartition $(X_1, Y_1)$
such that 
$|X_1| = 2m$,
$|Y_1| = (r-2)m$,
all vertices in $X_1$ have degree $r-2$,
and 
all vertices in $Y_1$ have degree $2$.
\item
Let $H_2$ be an $(r-4)$-edge-connected 
essentially $(r-3)$-edge-connected 
bipartite graph 
with bipartition $(X_2, Y_2)$
such that 
$|X_2| = (r-4)m$,
$|Y_2| = (r-2)m$,
all vertices in $X_2$ have degree $r-2$,
and 
all vertices in $Y_2$ have degree $r-4$.
\end{itemize}

\begin{figure}[htbp]
\begin{center}
\includegraphics[scale=0.9]{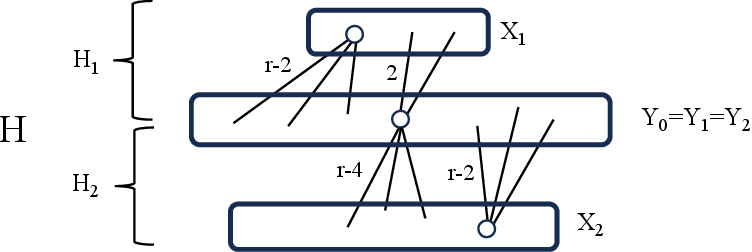}
\caption{An  $(r-2)$-regular bipartite graph $H$, where $|X_1|=2m, |Y_0|=|Y_1|=|Y_2|=(r-2)m, |X_2|=(r-4)m$.}
\label{fig-3A}
\end{center}
\end{figure}

Let $H$ be the bipartite graph 
obtained from $H_1$ and $H_2$ 
by identifying each vertex in $Y_1$ with a vertex in $Y_2$
along a bijection between $Y_1$ and $Y_2$.
We denote by $Y_0$ the set of vertices in $H$ obtained from $Y_1=Y_2$.
Note that $H$ is an $(r-2)$-regular bipartite graph
with bipartition $(X_1 \cup X_2, Y_0)$ such that $|X_1 \cup X_2| = |Y_0| = (r-2)m$
(see Fig.~\ref{fig-3A}).

\begin{figure}[htbp]
\begin{center}
\includegraphics[scale=0.9]{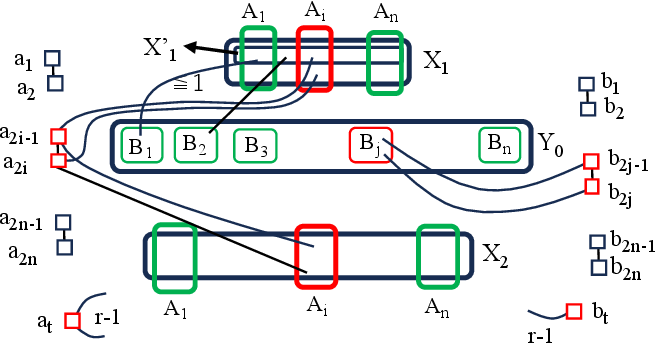}
\caption{A $K_{1,r}$-free $r$-edge connected $r$-regular graph, where
$(r-2)m=(r-1)n, |X_1'|=2n, |A_1|= \cdots  = |A_n|=r-1$, $|B_1|= \cdots = |B_n|=r-1$.}
\label{fig-3B}
\end{center}
\end{figure}

Let $n = (r-2)m/(r-1)$ be an integer.
Then $n < m$, and  we can take a subset $X_1' \subseteq X_1$ with $|X_1'| = 2n$.
Moreover, $|X_1\cup X_2|=(r-2)m =(r-1)n$, and thus we can
partition $X_1\cup X_2$ into $n$ sets $A_1, A_2, \dots , A_{n}$ 
with $|A_i \cap X_1'| = 2$ and $|A_i| = r-1$ for every $1 \leq i \leq n$.
Take $2n$ new vertices $a_1, a_2, \dots , a_{2n}$, and 
for each $1\le i \le n$,
we connect $a_{2i-1}$ and $a_{2i}$
to all vertices in $A_i$,
and further add an edge $a_{2i-1}a_{2i}$.
Then, in the resulting graph, every vertex in $A_1 \cup \cdots \cup A_{n}\cup \{a_1, \ldots, a_{2n}\}$ has degree $r$  (see Fig.~\ref{fig-3B}).

It follows from $n = (r-2)m/(r-1)$ and $m \geq 2(r-1)^2$ that 
\begin{align}
|X_1 - X_1'| = 2m - 2n
= \frac{2m}{r-1} \geq 4(r-1).
\label{eq-10}
\end{align}
Since all vertices in $Y_1$ have degree $2$ in $H_1$, it follows from (\ref{eq-10}) that
there are at least $2(r-1)$ vertices $y$ of $Y_0$ which are adjacent to $X_1-X_1'$ in $H_1$, and these vertices $y$ are adjacent to at most one vertex of $X_1'$ in $H_1$.
Let $B_1$ and $B_2$ be two disjoint sets of such vertices $y$
with $|B_1| = |B_2| = r-1$.
Then
\[ |Y_0 -(B_1 \cup B_2)|= (r-2)m -2(r-1)=(n-2)(r-1).\]
Thus, we can partition $Y_0 - (B_1 \cup B_2)$ into $n-2$ sets $B_3, B_4, \dots ,B_{n}$
with $|B_i| = r-1$ for every $3 \leq i \leq n$.
Take $2n$ new vertices $b_1, b_2, \dots , b_{2n}$,
and for each $1\le i \le n$,
we connect $b_{2i-1}$ and $b_{2i}$
to all vertices in $B_i$, and we further add an edge $b_{2i-1}b_{2i}$  (see Fig.~\ref{fig-3B}).

Let $G$ be the resulting graph. 
Since every vertex $x$ in $A_i$ belongs to the triangle $x a_{2i-1}a_{2i}$
and every vertex $y$ in $B_i$ belongs to the triangle $y b_{2i-1}b_{2i}$,
we see that $G$ is a $K_{1,r}$-free $r$-regular graph.

We now claim that $G$ is $r$-edge-connected.
Let $L$ be a minimal edge-cut of $G$, and 
let $D_1$ and $D_2$ be the two components of $G - L$.
If $D_1$ contains no vertex in $Y_0$,
then either for some $i$, $D_1$ consists of 
some vertices in $A_i$ possibly together with one or two of the vertices $a_{2i-1}$ and $a_{2i}$ (in the case where $V(D_1)\cap \{a_{2i-1}, a_{2i}\} \ne \emptyset$, it is possible that $V(D_1) \cap A_i=\emptyset$),
or one or two vertices of $b_j$'s.
In either case,
it is easy to see that $|L| \geq r$ by construction of $G$.
Thus, we may assume that $D_1$ contains a vertex in $Y_0$.
By symmetry,  we may also assume that $D_2$ also contains a vertex in $Y_0$.

Suppose that $V(D_1) \cap X_1\ne \emptyset$, $V(D_1) \cap X_2\ne \emptyset$,
$V(D_2) \cap X_1\ne \emptyset$ and $V(D_2) \cap X_2\ne \emptyset$.
Then for $i=1,2$,
$L \cap E(H_i)$ is an essential edge-cut of $H_i$,
and hence 
\[ |L| \geq |L \cap E(H_1)| + |L \cap E(H_2)|
\geq 3 + (r-3) = r,
\] 
where the inequality follows from the fact that
$H_1$ is essentially $3$-edge-connected 
and $H_2$ is essentially $(r-3)$-edge-connected.
We may therefore assume that at least one of $V(D_1) \cap X_1$, $V(D_1) \cap X_2$,
$V(D_2) \cap X_1$ and $V(G_2) \cap X_2$ is empty.
We here prove only the case $V(D_1) \cap X_1 = \emptyset$,
but the other cases can similarly be shown.
Assume that $V(D_1) \cap X_1 = \emptyset$. Then $X_1 \subset V(D_2)$, and 
every edge of $H_1$ incident with a vertex in $V(D_1) \cap Y_0$ belong to $L$.
Since $L \cap E(H_2)$ is an edge-cut of $H_2$,
we have $|L \cap E(H_2)| \geq r-4$.
If $|V(D_1) \cap Y_0| \geq 2$,
then 
we have 
\[ |L| \geq |L \cap E(H_1)| + |L \cap E(H_2)|
\geq 
2 |V(D_1) \cap Y_0| + (r-4)
\geq r.
\] 
Thus,
we may assume $|V(D_1) \cap Y_0| =1$.
The unique vertex in $V(D_1) \cap Y_0$ 
is adjacent to vertices $b_{2j-1}$ and $b_{2j}$
for some $1 \leq j \leq n$,
and  for each $h\in \{2j-1,2j\}$, at least one edge joining $b_h$ and $Y_0$ must  belong to $L$ since $|V(D_1) \cap Y_0| =1$.
We can now easily see that $|L| \geq r$.
Therefore,
$G$ is $r$-edge-connected as claimed.

Let $W = X_1' \cup \{b_1, b_3\}$.
Since $|X_1'| = 2n$, we see $|W|$ is even.
By the choice of $H_1, B_1$ and $B_2$,
we see that 
for every vertex $v$ of $G$, 
$N_G(v)$ contains at most two vertices of $W$.
Define a function $f :V(G) \to \mathbb{Z}^+$ as 
\[ f(v)= \bigg\{
\begin{array}{ll} 
1 & \mbox{if $v \in W$,~~ and } \\ 
2 & \mbox{otherwise.}
\end{array} 
    \] 
We let 
\begin{align*}
S  = & X_1 \cup X_2 \cup \{b_i : 1 \leq i \leq 2n\}, \quad \mbox{and}\\
T  = & Y_0 \cup \{a_i : 1 \leq i \leq 2n\}.
\end{align*} 
Note that $|S| = |X_1 \cup X_2| + 2n = |Y_0| + 2n= |T|$
and $|X_1'| = 2n$.
This implies 
\begin{align*}
\delta(S,T) & = f(S) + \deg_{G-S}(T) -f(T)  -q(S,T) \\
& = 2|S| - (|X_1'|+2) + 2n - 2|T| - 0  \\
& = -2,
\end{align*}
and hence it follows from Theorem \ref{thm-10} that 
$G$ has no spanning path-cycle system with respect to $W$.

\medskip
We next consider the case of $r=4$. Let $n\ge 6$ be an integer, and let $H=(x_1, y_1, x_2, y_2, \ldots , x_{3n}, y_{3n} , x_1)$ be a cycle of order $6n$, and let $X_1'=\{x_i: n+1 \le i \le 3n\}$. Let $A_i=\{x_i, x_{i+n}, x_{i+2n}\}$ and
$B_i=\{y_{3i-2}, y_{3i-1}, y_{3i}\}$ for every $1\le i \le n$. We construct a 4-edge-connected 
4-regular graph $G$ from the cycle $H$ and $4n$ new vertices $a_1,a_2, \ldots, a_{2n}, b_1, b_2, \ldots, b_{2n}$ by adding new edges as follows: Add edges $a_{2j-1}a_{2j}$ and $b_{2j-1}b_{2j}$ for every $1\le j \le n$, and for every $1 \le i \le n$, add edges between $\{a_{2i-1}, a_{2i}\}$ and $\{x_i, x_{n+i}, x_{2n+i}\}$, and between  $\{b_{2i-1}, b_{2i}\}$ and $\{y_{3i-2}, y_{3i-1}, y_{3i}\}$ (see Fig.~\ref{fig-6}). 
Let $W=\{x_{n+1},x_{n+2}, \ldots, x_{3n}\}\cup \{b_1,b_3\}$, and let $S=\{x_1,x_2, \ldots, x_{3n}\} \cup \{b_1,b_2, \ldots, b_{2n}\}$ and $T=\{y_1,y_2, \ldots, y_{3n}\}\cup \{a_1,a_2, \ldots, a_{2n}\}$. Define $f:V(G)\to \{1,2\}$ as $f(v)=1$  for $v\in W$ and otherwise $f(v)=2$. Then $|S|=|T|=5n$ and $|W|=2n+2$, and we have
\begin{align*}
\delta(S,T) & = f(S) + \deg_{G-S}(T) -f(T) -q(S,T) \\
 & = 2|S| -|W| + 2n -2|T| -0  \\
 & = -2<0.
 \end{align*}
Hence $G$ has no spanning path-cycle system with respect to $W$.

\begin{figure}[htbp]
\begin{center}
\includegraphics[scale=0.9]{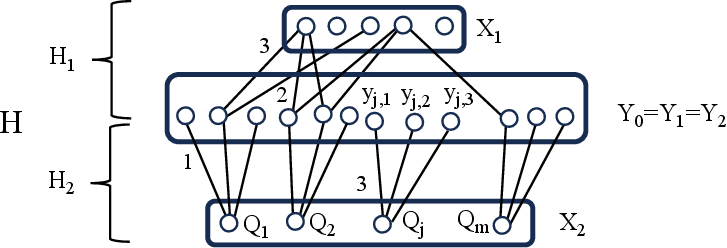}
\caption{A 3-regular bipartite graph $H$, in which $|X_1|=2m$, $|Y_0|=|Y_1|=|Y_2|=3m$ and $|X_2|=m$. }
\label{fig-6A}
\end{center}
\end{figure}

\begin{figure}[htbp]
\begin{center}
\includegraphics[scale=0.9]{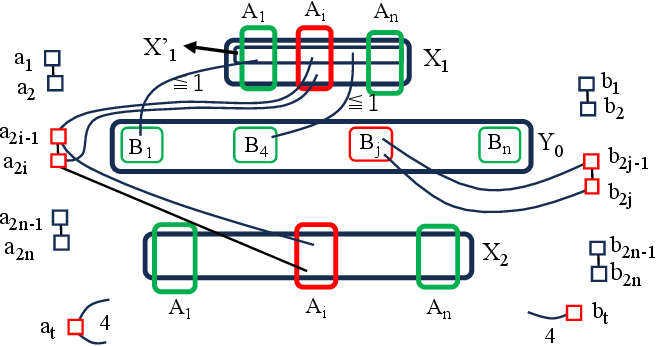}
\caption{A $K_{1,5}$-free 5-edge-connected 5-regular graph $G$, in which $3m=2n$, $|X'_1|=2n$, 
$|A_1|= \cdots =|A_n|=4$ and $|B_1|= \cdots =|B_n|=4$. }
\label{fig-6B}
\end{center}
\end{figure}

\medskip
Finally, we consider the case of $r=5$. In this case, let $m$ be a multiple of $(r-1)=4$ such that $m\ge 6\cdot 4^2$. As in the case where $r\ge 6$,  we construct a $K_{1,5}$-free $5$-regular $G$ as follows (see Figs.~\ref{fig-6A} and \ref{fig-6B}).  
Let $H_1$ be an essentially 3-edge-connected bipartite graph with bipartition $(X_1,Y_1)$ such that $|X_1|=2m$ and $|Y_1|=3m$, all vertices in $X_1$ have degree 3, and all vertices in $Y_1$ have degree 2, and let $H_2$ be a bipartite graph with bipartition $(X_2,Y_2)$ such that $|X_2|=m$ and $|Y_2|=3m$, all vertices in $X_2$ have degree 3, and all vertices in $Y_2$ have degree 1. Then $H_2$ consists of $m$ disjoint copies of $K_{1,3}$. Define $H$ and $Y_0$ as in the case where $r\ge 6$. 
Let $Q_1, \ldots, Q_m$ be the components of $H_2$, and for each $1\le j \le m$, write $V(Q_j)\cap Y_0=\{y_{j,1} , y_{j,2} , y_{j,3}\}$.

Let $n = (r-2)m/(r-1) = 3m/4$. Since $n < m$, we can take a subset $X_1' \subseteq X_1$ with $|X_1'| = 2n$.
By $m \geq 6\cdot 4^2$, we have
\begin{align}
|X_1 - X_1'| = 2m - 2n
= \frac{m}{2} \geq \frac{6\cdot 4^2}{2}=48.
\label{eq-16}
\end{align}

\begin{figure}[htbp]
\begin{center}
\includegraphics[scale=0.9]{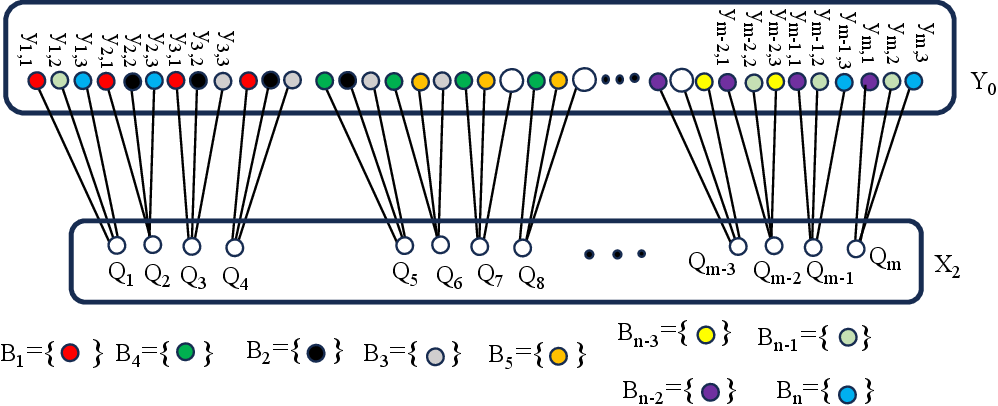}
\caption{A partition $B_1, B_2, B_3, B_4, B_5, \cdots, B_{n}$ of $Y_0$.  }
\label{fig-9}
\end{center}
\end{figure}

 As in the case where $r\ge 6$, we can take two disjoint subsets $B_1,B_4 \subset Y_0$ with $|B_1|=|B_4|=4$ so that each vertex in $B_1\cup B_4$ is adjacent to at most one vertex of $X_1'$. Further, since there are at least $24$ such vertices in $Y_0$ by (\ref{eq-16}), we can choose $B_1$ and $B_4$ so that the 8 vertices in $B_1\cup B_4$ belong to distinct components of $H_2$. We may assume that $B_1=\{y_{j,1} : 1\le j \le 4\}$ and $B_4=\{y_{j,1} : 5 \le j \le 8\}$. 
For each $1\le i \le n$ with $i\ne 1,4$, write $i=3p+h, h\in \{1,2,3\}$, and set $B_i=\{ y_{j,h} : 4p+h \le j \le 4p+h+3\}$ (first indices of $y$ are to be read modulo $m$). 
For example,  $B_2=\{y_{2,2}, y_{3,2}, y_{4,2}, y_{5,2}\}$,  $B_3=\{ y_{3,3} , y_{4,3}, y_{5,3}, y_{6,3}\}$, $B_5=\{ y_{6,2}, y_{7,2}, y_{8,2}, y_{9,2} \}$, and $B_6=\{y_{7,3}, y_{8,3}, y_{9,3}, y_{10,3}\}$.
Since $n=\frac{3}{4}m$ is multiple of $3$,
we see $B_n=\{y_{m-1,3}, y_{m,3}, y_{1,3}, y_{2,3}\}$,  $B_{n-1}=\{ y_{m-2,2} , y_{m-1,2}, y_{m,2}, y_{1,2}\}$, $B_{n-2}=\{ y_{m-3,1}, y_{m-2,1}, y_{m-1,1}, y_{m,1} \}$, and $B_{n-3}=\{y_{m-5,3}, y_{m-4,3}, y_{m-3,3}, y_{m-2,3}\}$ (see Fig.~\ref{fig-9}).

 Then $B_1, B_2, \ldots, B_n$ form a partition of $Y_0$.  By applying the same construction as in the proof of the case where $r\ge 6$, we can obtain the desired $K_{1,5}$-free 5-regular graph $G$.

Arguing as in the case where $r\ge 6$, we can show that $G$ is 5-edge-connected, and also show that the neighborhood of every vertex of $G$ contains at most two vertices of $W$, and $G$ has no spanning path-cycle system with respect to $W=X_1' \cup \{b_1,b_4\}$.

\bigskip \noindent
{\bf Acknowledgments} The authors would like to thank Akira Saito for his valuable suggestions. 

\bigskip \noindent
{\bf Funding}   The third author was supported by JSPS KAKENHI Grant Numbers 22K19773 and 23K03195.  This work was supported by the Research Institute for Mathematical Sciences, an International Joint Usage/Research Center located in Kyoto University.

\bigskip \noindent
{\bf Data Availability} No data was used for the research described in the article.

\bigskip \noindent 
{\large \bf Declarations}

\bigskip
\noindent
{\bf Conflict of interest} The authors have no relevant financial or non-financial interests to disclose.


\begin{thebibliography}{9}  
\bibitem{AK-2011} Akiyama, J.,  Kano, M. : 
{\em Factors and Factorizations of Graphs}, 
{\bf  LNM 2031} (Springer), (2011). 


\bibitem{FK-2023}
Furuya, M., M. Kano, M. :
Factors with red-blue coloring of claw-free graphs and cubic graphs,
{\em Graphs Combin.}, {\bf 39} (2023) \#85.
 
\bibitem{EKO-2025} Egawa, Y., Kano, M., Ozeki, K., 
 Spanning path-cycle systems with given end-vertices in regular graphs
 (full version), arXive: ??? 
 
\bibitem{FK-2024}
Furuya, M., M. Kano, M. :
Degree factors with red-blue coloring of regular graphs,
{\em Electron. J. Combin.}, {\bf 31 (1) } (2024) \#P1.40.

\bibitem{Kaiser-2008} Kaiser, T. :
 Disjoint $T$-paths in tough graphs. 
 {\em J. Graph Theory}, {\bf 59}  (2008)  1–10.

\bibitem{LK-2020} Lu, H.,  Kano, M. :
 Characterization of 1-tough graphs using factors. 
 {\em Discrete Math.}, {\bf 343} (2020) 111901.

\bibitem{Tutte-1952} Tutte, W.T. :
The factors of graphs, 
{\it Can. Ann. J. Math.},  {\bf 4} (1952) 314--328.

\bibitem{Tutte-1954} Tutte, W.T. :
A short proof of the factor theorem for finite graphs,  
{\it Can. Ann. J. Math.},  {\bf 4} (1954) 347--352.

\end{thebibliography}
\end{document}